\input amstex
\input amsppt.sty   
\hsize 26pc
\vsize 42pc
\def\nmb#1#2{#2}         
\def\ign#1{}             

\redefine\o{\circ}

\define\al{\alpha}

\define\de{\delta}

\define\la{\lambda}

\define\om{\omega}

\define\row#1#2#3{#1_{#2},\ldots,#1_{#3}}
\define\x{\times}
\define\pr{\operatorname{pr}}
\def\today{\ifcase\month\or
 January\or February\or March\or April\or May\or June\or
 July\or August\or September\or October\or November\or December\fi
 \space\number\day, \number\year}
\topmatter
\title More smoothly real compact spaces  
\endtitle
\author Andreas Kriegl 
\\Peter W. Michor\endauthor
\leftheadtext{A\. Kriegl, P\. W\. Michor}
\affil
Institut f\"ur Mathematik, Universit\"at Wien,\\
Strudlhofgasse 4, A-1090 Wien, Austria.
\endaffil
\address
{Institut f\"ur Mathematik, Universit\"at Wien,
Strudlhofgasse 4, A-1090 Wien, Austria}
\endaddress
\email KRIEGL\@AWIRAP.BITNET, MICHOR\@AWIRAP.BITNET \endemail
\date {1991} \enddate
\abstract{A topological space $X$ is called $\Cal A$-real compact, if 
every algebra homomorphism from $\Cal A$ to the reals is an 
evaluation at some point of $X$, where $\Cal A$ is an algebra of 
continuous functions. Our main interest lies on algebras of smooth 
functions. In \cite{AdR} it was shown that any separable Banach space 
is smoothly real compact. Here we generalize this result to a huge 
class of locally convex spaces including arbitrary products of 
separable Fr\'echet spaces.}
\endabstract
\endtopmatter


\document

In \cite{KMS} the notion of real compactness was generalized, by 
defining a topological space $X$ to be $\Cal A$-real-compact, if 
every algebra homomorphism $\al:\Cal A\to\Bbb R$ is just the 
evaluation at some point $a\in X$, where $\Cal A$ is a some 
subalgebra of $C(X,\Bbb R)$. In case $\Cal A$ equals the algebra 
$C(X,\Bbb R)$ of all continuous functions this condition reduces to 
the usual real-compactness. Our main interest lies on algebras 
$\Cal A$ of smooth functions.
In particular we showed in \cite{KMS} that every space admitting 
$\Cal A$-partitions of unity is $\Cal A$-real-compact. Furthermore any 
product of the real line $\Bbb R$ is $C^\infty$-real-compact.
A question we could not solve was, whether $\ell^1$ is 
$C^\infty$-real-compact, despite the fact that there are no smooth 
bump functions. \cite{AdR} had already shown that this is true not only for 
$\ell^1$, but for any separable Banach space.

The aim of this paper is to generalize this result of \cite{AdR} 
to a huge class of locally convex spaces, 
including arbitrary products of separable Fr\'echet spaces.

\subheading{Convention} All subalgebras $\Cal A\subseteq C(X,\Bbb R)$ 
are assumed to be real algebras with unit and with the additional 
property that for any $f\in\Cal A$ with $f(x)\ne 0$ for all $x\in X$ 
the function $\frac1f$ lies also in $\Cal A$.

\proclaim{\nmb.{1}. Lemma} Let $\Cal A\subset C(X,\Bbb R)$ be a finitely 
generated subalgebra 
of 
continuous functions on a topological space $X$. 
Then $X$ is $\Cal A$-real-compact.
\endproclaim
\demo{Proof} Let 
$\al:\Cal A\to\Bbb R$ be an algebra homomorphism.
We first show that for any finite set $\Cal F\subset\Cal A$
there exists a point $x\in X$ with $f(x)=\al(f)$ for 
all $f\in \Cal F$.

For $f\in \Cal A$ let 
$Z(f):=\{x\in X:f(x)=\al(f)\}$. 
Then 
$Z(f)=Z(f-\al(f)1)$, since $\al(f-\al(f)1)=0$. Hence we may assume that 
all $f\in\Cal F$ are even contained in $\ker\al=\{f:\al(f)=0\}$. Then 
$\bigcap_{f\in\Cal F}Z(f) = Z(\sum_{f\in\Cal F}f^2)$.
The sets $Z(f)$ are not empty, since otherwise $f\in\ker\al$ and 
$f(x)\ne 0$ for all $x$, so $\frac1f\in\Cal A$ and hence 
$1=f\frac1f\in\ker\al$, a contradiction to $\al(1)=1$.

Now the lemma is valid, whether the condition ``finitely 
generated'' is meant in the sense of an ordinary algebra or even as 
an algebra 
with the additional assumption on non-vanishing functions, since then any 
$f\in\Cal A$ can be written as a rational function in the elements of 
$\Cal F$. Thus $\al$ applied to such a rational function is just the 
rational function in the corresponding elements of 
$\al(\Cal F)=\Cal F(x)$, and is thus the value of the rational 
function at $x$.
\qed\enddemo

\proclaim{\nmb.{2}. Corollary} Any  
algebra-homomorphism $\al:\Cal A\to\Bbb R$ is monotone.
\endproclaim
\demo{Proof}
Let $f_1\leq f_2$. By \nmb!{1} there exists an $x\in X$ such that 
$\al(f_i)=f_i(x)$ for $i=1,2$. Thus 
$\al(f_1)=f_1(x)\leq f_2(x)=\al(f_2)$.
\qed\enddemo

\proclaim{\nmb.{3}. Corollary} Any algebra-homomorphism 
$\al:\Cal A\to\Bbb R$ is bounded, for every convenient algebra structure on 
$\Cal A$.
\endproclaim
By a convenient algebra structure we mean a convenient vector space 
structure for which the multiplication $\Cal A\x\Cal A\to\Cal A$ a bilinear 
bornological mapping. A convenient vector space is a separated locally 
convex vector space which is Mackey complete, see \cite{FK}.

\demo{Proof}
Suppose that $f_n$ is a bounded sequence, but $|\al(f_n)|$ is unbounded. 
Replacing $f_n$ by $f_n^2$ we may assume that $f_n\geq 0$ and hence
also $\al(f_n)\geq 0$. Choosing a subsequence we 
may even assume that $\al(f_n)\geq 2^n$. Now consider 
$\sum_{n}\frac1{2^n}f_n$. This series converges in the sense of Mackey, 
and since the bornology on $\Cal A$ is complete the limit is an element 
$f\in\Cal A$. Applying $\al$ yields
$$
\align
\al(f) &= \al\Biggl(\sum_{n=0}^N \frac1{2^n}f_n + 
          \sum_{n>N}\frac1{2^n}f_n\Biggr) 
     = \sum_{n=0}^N \frac1{2^n}\al(f_n) + 
          \al\Biggl(\sum_{n>N}\frac1{2^n}f_n\Biggr) \geq \\
     &\geq \sum_{n=0}^N \frac1{2^n}\al(f_n) + 0 
     = \sum_{n=0}^N \frac1{2^n}\al(f_n),
\endalign
$$
where we applied to the function
$\sum_{n>N}\frac1{2^n}f_n\geq 0$ that $\al$ is 
monotone. Thus the series $\sum_{n=0}^N \frac1{2^n}\al(f_n)$ is 
bounded and increasing, hence converges, but its summands are bounded 
by 1 from below. This is a contradiction.
\qed\enddemo

\subheading{\nmb.{4}. Definition}
We recall that a mapping $f:E\to F$ 
between convenient vector spaces is called smooth
($C^\infty$ for short), if the composite $f\o c:\Bbb R\to F$ is smooth  
for every smooth curve $c:\Bbb R\to E$. It can be shown that under these 
assumptions derivatives $f^{(p)}:E\to L^p(E,F)$ exist. See \cite{FK}.

A mapping is called $C^\infty_c$, if in addition all derivatives considered 
as mappings $d^pf:E\x E^p\to F$ are continuous.

\smallpagebreak
Now we generalize Lemma 5 and Proposition 7 of \cite{AdR} to 
arbitrary convenient vector spaces.

\subheading{\nmb.{5}. Definition} Let $\Cal A\subseteq C(X,\Bbb R)$ be 
a set of continuous functions on $X$. We say say that a space $X$ 
admits large carriers of class $\Cal A$, if for every neighborhood 
$U$ of a point $p\in X$ there exists a function $f\in\Cal A$ with 
$f(p)=0$ and $f(x)\ne 0$ for all $x\notin U$.

Every $\Cal A$-regular space $X$ admits large $\Cal A$-carriers, where 
$X$ is called $\Cal A$-regular if for every neighborhood $U$ of a 
point $p\in X$ there exists a function $f\in\Cal A$ with $f(p)>0$ and 
$f(x)=0$ for $x\notin U$. The existence of large $\Cal A$-carriers 
follows by using the modified function $\bar f:=f(a)-f$.

In \cite{AdR, Proof of theorem 8} it is proved, that every separable 
Banach space admits large $C^\infty_c$-carriers.
The carrying functions can even be chosen as polynomials as shown in 
lemma \nmb!{7} below.

\proclaim{\nmb.{6}. Lemma} Let $E$ be a convenient vector space, 
$\{x_n':n\in\Bbb N\}\subset E'$ be bounded, $(\la_n)\in\ell^1(\Bbb N)$ 
Then the series 
$(x,y)\mapsto\sum_{n=1}^\infty\la_n x_n'(x)x'_n(y)$ converges to a continuous 
symmetric bilinear function on $E\x E$.
\endproclaim

\demo{Proof} Clearly the function converges pointwise.
Since the sequence $\{x_n'\}$ is bounded, it is 
equicontinuous, hence bounded on some neighborhood $U$ of 0, so there 
exists a constant $M\in\Bbb R$ such that $|x_n'(U)|\leq M$ for all 
$n\in \Bbb N$. For $x,y\in U$ we have 
$|\sum_{n=1}^\infty\la_n x_n'(x)x'_n(y)|\le \sum_{n=1}^\infty|\la_n|M^2$, 
which suffices for continuity of a bilinear function.
\qed\enddemo

\proclaim{\nmb.{7}. Lemma} Let $E$ be a Banach space which is 
separable or whose dual is separable for the topology of pointwise 
convergence. Then $E$ admits large 
carriers for continuous polynomials  of degree 2.
\endproclaim

\demo{Proof}
If $E$ is separable
there exists a dense sequence $(x_n)$ in $E$.
By the Hahn-Banach theorem 
\cite{J, 7.2.4} there exist $x_n'\in E'$ with 
$x_n'(x_n)=|x_n|$ and $|x_n'|\leq 1$. 

Claim: $\sup_n |x_n'(x)|=|x|$ \newline
Since $|x_n'|\leq 1$ we have $(\leq)$. For the converse direction 
let $\de>0$ be given. By denseness there exists an $n\in\Bbb N$ such that 
$|x_n-x|<\tfrac{\de}2$. So we have:
$$
\align
|x| &\leq |x_n| + |x-x_n| < |x_n'(x_n)| + \tfrac\de2 \leq \\
     &\leq |x_n'(x)| + \undersetbrace < |x-x_n| < \tfrac\de2 
   	\to {|x_n'(x-x_n)|} + \tfrac\de 2 < \\
     &< |x_n'(x)| + \de.
\endalign
$$
If the dual $E'$ is separable for the topology of pointwise 
convergence, then let $x_n'$ be a sequence which is weakly dense in 
the unit ball of $E'$. Then $|x|=\sup_n |x_n'(x)|$.

In both cases the continuous polynomials of lemma \nmb!{6} 
$$x\mapsto\sum_{n=1}^\infty\frac1{n^2} x_n'(x-a)^2$$ 
vanish exactly at $a$.
\qed\enddemo

\proclaim{\nmb.{8}. Lemma} Let $\al:\Cal A\to \Bbb R$ be an algebra 
homomorphism and assume that some subset $\Cal A_0\subset\Cal A$ 
exists and a point $a\in X$ such that $\al(f_0)=f_0(a)$ for all 
$f_0\in \Cal A_0$ and such that $X$ admits large carriers of class $\Cal A_0$.

Then $\al(f)=f(a)$ for all $f\in\Cal A$.
\endproclaim
\demo{Proof}
Let $f\in\Cal A$ be arbitrary. Since 
$X$ admits large $\Cal A_0$-carriers there exists 
for every neighborhood $U$ of $a$ 
a function $f_U\in\Cal A_0$ with $f_U(a)=0$ and $f_U(x)\ne 0$ for all 
$x\in U$. By lemma \nmb!{1} there exists a point 
$a_U$ such that $\al(f)=f(a_U)$ and $\al(f_U)=f_U(a_U)$. Since 
$f_U\in \Cal A_0$, we have $f_U(a_U)=\al(f_U)=f_U(a)=0$, hence 
$a_U\in U$. Thus the net $a_U$ converges to $a$ and consequently  
$f(a)=f(\lim_Ua_U)=\lim_U f(a_U)=\lim_U \al(f)=\al(f)$ since $f$ is 
continuous.
\qed\enddemo

Now we generalize proposition 2 and lemma 3 of \cite{BBL}. 
Let for every convenient vector space $E$ a subalgebra $\Cal A(E)$ of 
$C(E,\Bbb R)$ be given, such that for every $f\in L(E,F)$ the image of $f^*$ 
on $\Cal A(F)$ lies in $\Cal A(E)$. Examples are $C^\infty_c$, 
$C^\infty\cap C$, $C^\om_c:=C^\infty_c\cap C^\om$, $C^\om\cap C$, 
where $C^\om$ denotes the algebra of real analytic functions in the 
sense of \cite{KM}, and suitable algebras of functions of inite 
differentiability like $\operatorname{Lip}^m$ (see \cite{FK}) or 
$C^m_c$.

\proclaim{\nmb.{9}. Theorem} Let $E_i$ be $\Cal A$-real-compact 
spaces that admit large carriers of class $\Cal A$. Then any closed subspace 
of the product of the spaces $E_i$, and in particular every projective limit 
of these spaces, has the same properties.
\endproclaim

\demo{Proof}
First we show that this is true for the product $E$. We use lemma \nmb!{8} 
with $\Cal A(E)$ for $\Cal A$ and the vector space generated by
$\bigcup_i \{f\o \pr_i:f\in \Cal A(E_i)\}$ for $\Cal A_0$, where 
$\pr_j:E=\prod_i E_i\to E_j$ denotes the canonical projection. 
Let the finite sum
$f=\sum_i f_i\o \pr_i$ be an element of $\Cal A_0$. Since  
$\al\o \pr_i{}^*:\Cal A(E_i)\to \Cal A(E)\to\Bbb R$ is an 
algebra homomorphism, there exists a point $a_i\in E_i$ such 
that $\al(f_i\o\pr_i)=(\al\o \pr_i{}^*)(f_i)=f_i(a_i)$. 
Let $a$ be the point in $E$ with 
coordinates $a_i$. Then
$$
\align
\al(f)&=\al(\sum_i f_i\o\pr_i)
     =\sum_i\al(f_i\o\pr_i)\\
     &=\sum_i f_i(a_i)
     =\sum_i (f_i\o\pr_i)(a)
     =f(a)
\endalign
$$
Now let $U$ be a neighborhood of $a$ in $E$. Since we consider the 
product topology on $E$ we may assume that 
$a\in\prod U_i\subset U$, where $U_i$ are neighborhoods of $a_i$ 
in $E_i$ and are equal to $E_i$ except for $i$ in some finite subset 
$F$ of the index set.
Now choose $f_i\in \Cal A(E_i)$ with $f_i(a_i)=0$ and 
$f_i(x_i)\ne 0$ for all $x_i\notin U_i$. Consider 
$f=\sum_{i\in F}(f_i\o\pr_i)^2\in \Cal A_0$. Then 
$f(a)=\sum_{i\in F}f_i(a_i)^2=0$. Furthermore $x\notin U$ implies that 
$x_i\notin U_i$ for some $i$, which turns out to be in $F$, and hence 
$f(x)\geq f_i(x_i)^2>0$. So we may apply lemma \nmb!{8} to conclude that
$\al(f)=f(a)$ for all $f\in \Cal A(E)$.

Now we prove the result for a closed subspace $F\subset E$. Again we 
want to apply lemma \nmb!{8}, this time with 
$\Cal A(F)$ for $\Cal A$ and $\{f|_F:f\in \Cal A(E)\}$ 
for $\Cal A_0$. Since 
$\al\o\operatorname{incl}^*:\Cal A(E)\to \Cal A(F)\to\Bbb R$ is an algebra 
homomorphism there exists an $a\in E$ with $\al(f|_F)=f(a)$ for all 
$f\in \Cal A(E)$. Now let $U$ be a neighborhood of $a$ in $E$ then 
there exists an $f_U\in \Cal A(E)$ with $f_U(a)=0$ and 
$f_U(x)\ne 0$ for all $x\notin U$. 
By lemma \nmb!{1} there exists a point $a_U\in F$ such that 
$f_U(a_U)=\al(f_U|_F)=f_U(a)=0$. Hence $a_U$ is 
in $U$, and thus is a net in $F$ which converges to $a$. In 
particular $a\in F$, since $F$ is 
closed in $E$.
If $V$ is a neighborhood of $a$ in $F$ then there exists a 
neighborhood $U$ of $a$ in $E$ with $U\cap F\subset V$ and hence an 
$f\in\Cal A_0$ with $f(a)=0$ and $f(x)\ne 0$ for all $x\notin U$. So 
again \nmb!{8} applies.
\qed\enddemo

\subheading{\nmb.{10}. Remark}
Theorem \nmb!{9} shows that a closed subspace of a product
of certain $\Cal A$-real-compact spaces is again $\Cal A$-real-compact. 
Of course the natural question arises, whether the result remains 
true for arbitrary $\Cal A$-real-compact spaces.

It is even open, whether the product of two $\Cal A$-real-compact spaces is 
$\Cal A$-real-compact, or whether a closed subspace of an 
$\Cal A$-real-compact space is $\Cal A$-real-compact, or whether
a projective limit of a projective 
system of $\Cal A$-real-compact spaces is $\Cal A$-real-compact.

\proclaim{\nmb.{11}. Corollary} Let $E$ be a separable Fr\'echet space (e.g.
a Fr\'echet-Montel space), then 
every algebra homomorphism on $C^\infty(E,\Bbb R)$ or on $C^\infty_c(E,\Bbb R)$ is a 
point evaluation. The same is true for any product of separable 
Fr\'echet spaces.
\endproclaim

\demo{Proof}
Any Fr\'echet space has a countable Basis $\Cal U$ of absolutely 
convex 0-neigh\-bor\-hoods, 
and since it is complete it is a closed subspace of the product 
$\prod_{u\in\Cal U} \widetilde{E_{(U)}}$. The $E_{(U)}$ are the normed 
spaces formed by $E$ modulo the kernel of the Minkowski functional 
generated by $U$. As quotients of 
$E$ the spaces $E_{(U)}$ are separable if $E$ is such. So the completion 
$\widetilde{E_{(U)}}$ is a separable Banach space and hence by 
\cite{AdR, Theorem 8} $\widetilde{E_{(U)}}$ is 
$C^\infty_c$-real-compact and admits large 
$C^\infty_c$-carriers. By theorem \nmb!{9} the same is true for the given 
Fr\'echet space.
So the result is true for $C^\infty_c(E,\Bbb R)$. 
Since $E$ is metrizable this algebra coincides with $C^\infty(E,\Bbb R)$, see 
\cite{K, 82}.

Now for a product $E$ of metrizable spaces the two algebras 
$C^\infty(E,\Bbb R)$ and $C^\infty_c(E,\Bbb R)$ again coincide. This 
can be seen as follows. For every countable subset $A$ of the index 
set, the corresponding product is separable and metrizable, hence 
$C^\infty$-real-compact. Thus there exists a point $x_A$ in this 
countable product such that $\al(f)=f(x_A)$ for all $f$ which factor 
over the projection to that countable subproduct. Since for 
$A_1\subset A_2$ the projection of $x_{A_2}$ to the product over 
$A_1$ is just $x_{A_1}$ (use the coordinate projections composed with 
functions on the factors for $f$), there is a point $x$ in the 
product, whose projection to the subproduct with index set $A$ is 
just $x_A$. Every Mackey continuous function, and in particular every 
$C^\infty$-function, depends only on countable many coordinates, thus 
factors over the projection to some subproduct with countable index 
set $A$, hence $\al(f)=f(x_A)=f(x)$. This can be shown by the same 
proof as for a product of factors $\Bbb R$ in \cite{FK, Theorem 
6.2.9}, since the result of \cite{M, 1952} is valid for a product of 
separable metrizable spaces.
\qed\enddemo

\enddocument